\documentclass[leqno,a4paper]{article}

\usepackage{amssymb,amsmath,amsthm}

\newtheorem{theorem}{Theorem}

\newtheorem{lemma}[theorem]{Lemma}

\newtheorem{remark}[theorem]{Remark}

\begin{document}
\title{\textbf{First-order invariants}\\\textbf{of differential }$2$\textbf{-forms}}
\author{\textsc{J. Mu\~{n}oz Masqu\'{e}}\thanks{Instituto de Tecnolog\'{\i}as F\'{\i}sicas y de la
Informaci\'on, Consejo Superior de Investigaciones Cient\'{\i}ficas,
 C/ Serrano 144, 28006-Madrid (Spain), email: jm.masque@gmail.com}
\and \textsc{L. M. Pozo Coronado}\thanks{Departamento de Matem\'atica Aplicada a las TIC,
ETSI Sistemas Inform\'aticos, Universidad Polit\'ecnica de Madrid,
C\ Alan Turing, s/n, 28031-Madrid (Spain), email: lm.pozo@upm.es}}

\date{}
\maketitle

\begin{abstract}
\noindent Let $M$ be a smooth manifold of dimension $2n$, and let $O_{M}$ be
the dense open subbundle in $\wedge^{2}T^{\ast}M$ of $2$-covectors of maximal
rank. The algebra of $\operatorname*{Diff}M$-invariant smooth functions of
first order on $O_{M}$ is proved to be isomorphic to the algebra of smooth
$Sp(\Omega_{x})$-invariant functions on $\wedge^{3}T_{x}^{\ast}M$, $x$ being a
fixed point in $M$, and $\Omega_{x}$ a fixed element in $(O_{M})_{x}$. The
maximum number of functionally independent invariants is computed.
\end{abstract}

\bigskip

\noindent\emph{Mathematics Subject Classification 2010:\/} Primary: 53A55
; Secondary: 22E15,
53D05,
58A10,
58A20

\smallskip

\noindent\emph{Key words and phrases:}\/Differential invariant function,
differential $2$-form, jet bundle, linear representation, symplectic group.

\bigskip

\section{Reduction to symplectic group\label{reduction}}
Let $M$ be a $C^{\infty}$ manifold and let
\[
G_{x}^{r}=G_{x}^{r}(M)=\left\{ j_{x}^{r}\phi\in J^{r}(M,M):\phi(x)=x,\det
\phi_{\ast,x}\neq0\right\} ,
\quad r\geq 1,
\]
be the Lie group of $r$-jets of diffeomorphisms at $x\in M$. If $r\geq s$, then
$G_{x}^{rs}$ denotes the kernel of the natural projection $G_{x}^{r}\to
G_{x}^{s}$.

In particular, for every $r\geq2$, $G_{x}^{r,r-1}$ is isomorphic to the
vectorial group $S^{r}T_{x}^{\ast}(M)\otimes T_{x}(M)$, as $J^{r}(M,M)\to
J^{r-1}(M,M)$ is an affine bundle modelled over
$S^{r}T^{\ast}(M)\otimes_{J^{r-1}(M,M)}T(M)$ and, therefore for every
$j_{x}^{r}\phi\in G_{x}^{r,r-1}$ there exists a unique $t\in
S^{r}T_{x}^{\ast}(M)\otimes T_{x}(M)$ such that
$t+j_{x}^{r}(1_{M})=j_{x}^{r}\phi$. Hence we can identify $j_{x}^{r}\phi$ to
$t$.
\begin{theorem}
\label{proposition_delta}
Let $M$ be a $C^{\infty}$ manifold of dimension $2n$, and let $p\colon O_{M}\to
M$ be the dense open subbundle in $\wedge^{2}T^{\ast}M$ of $2$-covectors of
maximal rank. Given a point $x\in M$, the map $\delta_{x}\colon
J_{x}^{1}O_{M}\to\wedge^{3}T_{x}^{\ast}M$,
$\delta_{x}(j_{x}^{1}\Omega)=(d\Omega)_{x}$, is a $G_{x}^{2}$-equivariant
$G_{x}^{21}$-invariant epimorphism.
\end{theorem}
\begin{proof}
The $G_{x}^{2}$-equivariance of $\delta_{x}$ is a consequence of the following
well-known property: $(d(\phi^{\ast}\Omega))_{x}=\phi^{\ast}((d\Omega
)_{\phi(x)})$. In fact, for every $\phi\in\operatorname*{Diff}_{x}M$ one
obtains
\begin{align*}
\delta_{x}\left(
j_{x}^{2}\phi\cdot j_{x}^{1}\Omega
\right)
&  =\delta _{x}
\left(
j_{x}^{1}\left( (\phi^{-1})^{\ast}\Omega
\right)
\right)
=\left( d((\phi^{-1})^{\ast}\Omega)\right) _{x}=(\phi^{-1})^{\ast}
\left(
(d\Omega)_{\phi^{-1}(x)}\right) \\ &
=j_{x}^{2}\phi\cdot\delta_{x}(j_{x}^{1}\Omega).
\end{align*}

Next, we show that $\delta_{x}$ is surjective. Let
$w_{3}\in\wedge^{3}T_{x}^{\ast}M$. Let $(x^{i})_{i=1}^{2n}$ be a coordinate
system centred at $x$. If
$w_{3}=\sum_{h<i<j}\lambda_{hij}(dx^{h})_{x}\wedge(dx^{i})_{x}
\wedge(dx^{j})_{x}$, then $w_{3}=\delta_{x}(j_{x}^{1}\Omega)$, where
$\Omega=\sum_{h<i<j}\lambda_{hij}x^{h}dx^{i}\wedge dx^{j}$. If
$(y_{hij})_{1\leq h<i<j\leq2n}$ is the coordinate system on
$\wedge^{3}T_{x}^{\ast}M$ given by
\[
w_{3}=\sum_{h<i<j}y_{hij}(w_{3})(dx^{h})_{x}\wedge(dx^{i})_{x}\wedge
(dx^{j})_{x},\quad\forall w_{3}\in\wedge^{3}T_{x}^{\ast}M,
\]
then the equations for $\delta_{x}$ are
$y_{abc}\circ\delta_{x}=y_{bc,a}-y_{ac,b}+y_{ab,c}$, $a<b<c$.

This proves that $\delta_{x}$ is the restriction of a linear mapping. Moreover,
the projection $p^{10}\colon J^{1}(\wedge^{2}T^{\ast}M)\to
J^{0}(\wedge^{2}T^{\ast}M)=\wedge^{2}T^{\ast}M$ is an affine bundle modelled
over $T^{\ast}M\otimes\wedge^{2}T^{\ast}M$, where the sum operation is defined
as follows:
\begin{equation}
(df)_{x}\otimes w_{2}+j_{x}^{1}\Omega=j_{x}^{1}\left[  \left(  f-f(x)\right)
\Omega^{\prime}+\Omega\right] , \label{operacion_suma1}
\end{equation}
$\Omega^{\prime}$ being any $2$-form such that
$\Omega_{x}^{\prime}=w_{2}\in\wedge_{x}^{2}T^{\ast}M$.

If $j_{x}^{2}\phi\in G_{x}^{21}$, then $(\phi^{\ast}\Omega)_{x}=\Omega_{x}$,
i.e, $p^{10}(j_{x}^{1}\Omega)=p^{10}(j_{x}^{1}(\phi^{\ast}\Omega))\,$; hence
there exists a unique $\tau=\sum_{j<k}\tau_{jkl}(dx^{l})_{x}\otimes
(dx^{j})_{x}\wedge(dx^{k})_{x}\in T_{x}^{\ast}M\otimes\wedge^{2}T_{x}^{\ast}M$
such that $j_{x}^{1}(\phi^{\ast}\Omega)=\tau+j_{x}^{1}\Omega$. If $\Omega
=\sum_{h<i}F_{hi}dx^{h}\wedge dx^{i}$, then
\[
\phi^{\ast}\Omega=\sum\nolimits_{j<k}\bar{F}_{jk}dx^{j}\wedge dx^{k},\quad
\bar{F}_{jk}=\sum\nolimits_{h<i}\left(  F_{hi}\circ\phi\right)  \det
\tfrac{\partial(\phi^{h},\phi^{i})}{\partial(x^{j},x^{k})},
\quad\phi^{h}=x^{h}\circ\phi,
\]
and taking derivatives for $\bar{F}_{jk}$ and evaluating at $x$, one obtains
\begin{align*}
\tfrac{\partial\bar{F}_{jk}}{\partial x^{l}}(x)
&  =\sum\nolimits_{h<i}
\tfrac{\partial F_{hi}}{\partial x^{a}}(\phi(x))
\tfrac{\partial\phi^{a}}{\partial x^{l}}(x)
\det\tfrac{\partial(\phi^{h},\phi^{i})}{\partial (x^{j},x^{k})}(x)\\
&  +\sum\nolimits_{h<i}F_{hi}(\phi(x))
\left\{
\tfrac{\partial^{2}\phi^{h}}{\partial x^{j}\partial x^{l}}(x)
\tfrac{\partial\phi^{i}}{\partial x^{k}}(x)
+\tfrac{\partial\phi^{h}}{\partial x^{j}}(x)
\tfrac{\partial^{2}\phi^{i}}{\partial x^{k}\partial x^{l}}(x)
\right. \\
&  \left. -\tfrac{\partial^{2}\phi^{h}}{\partial x^{k}\partial x^{l}}(x)
\tfrac{\partial\phi^{i}}{\partial x^{j}}(x)
-\tfrac{\partial\phi^{h}}{\partial x^{k}}(x)
\tfrac{\partial^{2}\phi^{i}}{\partial x^{j}\partial x^{l}
}(x)\right\} .
\end{align*}

As $j_{x}^{2}\phi\in G_{x}^{21}$, it follows: $\phi(x)=x$, $\frac{\partial
\phi^{h}}{\partial x^{i}}(x)=\delta_{i}^{h}$, and from the previous formula one
thus deduces $\tau_{jkl} =F_{hk}(x)\frac{\partial^{2}\phi^{h}}{\partial
x^{j}\partial x^{l}}(x) -F_{hj}(x)\frac{\partial^{2}\phi^{h}}{\partial
x^{k}\partial x^{l}}(x)$. As $\tau_{jkl}$ is alternate on $j,k$, one can write
$\tau=\frac{1}{2}\tau_{jkl}(dx^{l})_{x}\otimes(dx^{j})_{x}\wedge(dx^{k})_{x}$,
and recalling that the coordinates are centred at $x$, taking the formula
\eqref{operacion_suma1} into account, it follows:
$j_{x}^{1}\Omega^{\prime}=\tau$, where $\Omega^{\prime}$ is the $2$-form given
by
\begin{align}
\Omega^{\prime}  &  =\tfrac{1}{2}x^{l}\left\{  F_{hk}(x)\tfrac{\partial
^{2}\phi^{h}}{\partial x^{j}\partial x^{l}}(x)-F_{hj}(x)\tfrac{\partial
^{2}\phi^{h}}{\partial x^{k}\partial x^{l}}(x)\right\}  dx^{j}\wedge
dx^{k}\label{Omega_prima}\\ &  =\tfrac{1}{2}x^{l}d\left(
\tfrac{\partial\phi^{h}}{\partial x^{l}}\right)
\wedge\left\{  F_{hk}(x)dx^{k}+F_{hj}(x)dx^{j}\right\} \nonumber\\
&  =d\left\{  \tfrac{1}{2}\left(  x^{l}\tfrac{\partial\phi^{h}}{\partial
x^{l}}-\phi^{h}\right)  \left(  F_{hk}(x)dx^{k}+F_{hj}(x)dx^{j}\right)
\right\}  .\nonumber
\end{align}
Hence, $\delta_{x}\left(  j_{x}^{1}(\phi^{\ast}\Omega)\right)  =\delta
_{x}\left(  j_{x}^{1}\left(  \Omega^{\prime}+\Omega\right)  \right)
=\delta_{x}\left(  j_{x}^{1}\Omega\right)  $.
\end{proof}
Let $M$ be an arbitrary $C^{\infty}$-manifold and let $\bar{\phi}\colon
\wedge^{2}T^{\ast}M\to\wedge^{2}T^{\ast}M$ be the natural lift
of a diffeomorphism $\phi\in\operatorname*{Diff}M$; i.e.,
$\bar{\phi}(w)=(\phi^{-1})^{\ast}w$ for every $2$-covector
$w\in\wedge^{2}T^{\ast}M$. If $\Omega$ is a $2$-form on $M$, then
$\bar{\phi}\circ\Omega\circ\phi^{-1}=(\phi^{-1})^{\ast}\Omega$. Let
\[
\begin{array}
[c]{l} J^{1}\bar{\phi}\colon J^{1}\bigl(\wedge^{2}T^{\ast}M\bigr)\to
J^{1}\bigl(\wedge^{2}T^{\ast}M\bigr),\smallskip\\
J^{1}\bar{\phi}(j_{x}^{1}\Omega)=j_{\phi(x)}^{1}(\bar{\phi}\circ\Omega
\circ\phi^{-1}),
\end{array}
\]
be the $1$-jet prolongation of $\bar{\phi}$. A subset $S\subseteq
J^{1}(\wedge^{2}T^{\ast}M)$ is said to be natural if $\left( J^{1}\bar{\phi
}\right) (S)\subseteq S$ for every $\phi\in\operatorname*{Diff}M$. Let
$S\subseteq J^{1}(\wedge^{2}T^{\ast}M)$ be a natural embedded submanifold. A
smooth function $I\colon S\to\mathbb{R}$ is said to be an invariant of first
order under diffeomorphisms or even $\operatorname*{Diff}M$-invariant if
$I\circ J^{1}\bar{\phi}=I$, $\forall\phi\in\operatorname*{Diff}M$. If we set
$I(\Omega)=I\circ j^{1}\Omega$, for a given $2$-form $\Omega$ on $M$, then the
previous invariance condition reads as $I\left( (\phi
^{-1})^{\ast}\Omega)(\phi(x)\right) =I(\Omega)(x)$, for all $x\in
M$,$\;\phi\in\operatorname*{Diff}M$, thus leading us to the naive definition of
an invariant, as being a function depending on the coefficients of $\Omega$ and
its partial derivatives up to first order, which remains unchanged under
arbitrary changes of coordinates.
\begin{theorem}
\label{Th1} Let $M$ be a smooth connected manifold of dimension $2n$. The ring
of invariants of first order on $O_{M}$ is isomorphic to $C^{\infty}\left(
\wedge^{3}T_{x}^{\ast}M\right)  ^{Sp(\Omega_{x})}$, where $\Omega_{x}$ is a
fixed element in $O_{M}$.
\end{theorem}
\begin{proof}
As $M$ is connected, the group $\operatorname*{Diff}M$ acts transitively on
$M$. Therefore, it suffices to fix a point $x\in M$ and to compute
$\operatorname*{Diff}_{x}M$-invariant functions in
$C^{\infty}(J_{x}^{1}O_{M})$. From the very definitions it follows:
\[
C^{\infty}(J_{x}^{1}O_{M})^{\operatorname*{Diff}_{x}M}
=C^{\infty}(J_{x}^{1}O_{M})^{G_{x}^{2}},
\]
and by virtue of Theorem \ref{proposition_delta} the map
$\bar{\delta}_{x}\colon J_{x}^{1}O_{M}
\to(O_{M})_{x}\times\wedge^{3}T_{x}^{\ast}M$,
defined as follows:
$\bar{\delta}_{x}(j_{x}^{1}\Omega)=(\Omega_{x},(d\Omega)_{x})$, is
$G_{x}^{21}$-invariant and surjective; hence the induced homomorphism $\left(
\bar{\delta}_{x}\right)  ^{\ast}\colon C^{\infty}((O_{M})_{x}
\times\wedge^{3}T_{x}^{\ast}M)\to C^{\infty}(J_{x}^{1}O_{M})^{G_{x}^{21}}$
is injective. Next, we shall prove that $\left( \bar{\delta}_{x}\right)
^{\ast}$ is also surjective, by showing that every $I\in
C^{\infty}(J_{x}^{1}O_{M})^{G_{x}^{21}}$ takes constant value on the fibres of
$\bar{\delta}_{x}$, as in this case $I$ induces $\bar{I}\in
C^{\infty}((O_{M})_{x}\times\wedge^{3}T_{x}^{\ast}M)$ such that $I=\bar
{I}\circ\bar{\delta}_{x}=\left(  \bar{\delta}_{x}\right)  ^{\ast}(\bar{I})$.
Actually, this is a consequence of the fact that the fibres of $\bar{\delta
}_{x}$ coincide with the orbits of $G_{x}^{21}$ on $J_{x}^{1}O_{M}$. To prove
this, we first observe that every $j_{x}^{1}\bar{\Omega}\in\left(  \bar
{\delta}_{x}\right)  ^{-1}(\bar{\delta}_{x}(j_{x}^{1}\Omega))$ can be written
as $j_{x}^{1}\bar{\Omega}=j_{x}^{1}(\tilde{\Omega}+\Omega)$ with
$\tilde{\Omega}_{x}=0$, $(d\tilde{\Omega})_{x}=0$, and the proof reduces to
show the existence of a $2$-form $\Omega^{\prime}$ given by the formula
\eqref{Omega_prima} such that $j_{x}^{1}\Omega^{\prime}=j_{x}^{1}\tilde
{\Omega}$, since we have seen that
$j_{x}^{1}(\phi^{\ast}\Omega)=j_{x}^{1}\left(  \Omega^{\prime}+\Omega\right)
$, for some $j_{x}^{2}\phi\in G_{x}^{21}$. The rank of $\Omega$ being $2n$,
there exists a coordinate system $(x^{i})_{i=1}^{2n}$ centred at $x$ such that
$\Omega_{x}=\sum_{i=1}^{n}(dx^{2i-1})_{x}\wedge(dx^{2i})_{x}$, or equivalently,
$F_{2i-1,2i}(x)=1$, and $F_{jk}(x)=0$, $1\leq j<k\leq2n$ otherwise. We have
$\tilde{\Omega}=\sum\nolimits_{i<j}\lambda_{ijk}x^{k}dx^{i}\wedge dx^{j}+$
terms of order $\geq2$, because $\tilde{\Omega}$ vanishes at $x$, and by
imposing $\tilde{\Omega}$ to be closed at $x$, we obtain
\begin{equation}
0=\lambda_{ijk}-\lambda_{kji}+\lambda_{kij},\qquad1\leq k<i<j\leq2n.
\label{compat}
\end{equation}
Then, the equation $j_{x}^{1}\Omega^{\prime}=j_{x}^{1}\tilde{\Omega}$ is
equivalent to the system
\[
F_{hk}(x)\tfrac{\partial^{2}\phi^{h}}{\partial x^{j}\partial x^{l}}(x)
-F_{hj}(x)\tfrac{\partial^{2}\phi^{h}}{\partial x^{k}\partial x^{l}}(x)
=\lambda_{jkl,}\qquad1\leq j<k\leq2n,
\]
or equivalently, for $1\leq j^{\prime}<k^{\prime}\leq n$,
\begin{equation}
\begin{array}
[c]{rl}
\tfrac{\partial^{2}\phi^{2k^{\prime}-1}}
{\partial x^{2j^{\prime}}\partial x^{l}}(x)
-\tfrac{\partial^{2}\phi^{2j^{\prime}-1}}{\partial x^{2k^{\prime}}\partial
x^{l}}(x) = & \lambda_{2j^{\prime},2k^{\prime},l},
\smallskip \\
-\tfrac{\partial^{2}\phi^{2k^{\prime}}} {\partial x^{2j^{\prime}}\partial
x^{l}}(x) -\tfrac{\partial^{2}\phi^{2j^{\prime}-1}} {\partial
x^{2k^{\prime}-1}\partial x^{l}}(x) = & \lambda_{2j^{\prime},2k^{\prime}-1,l},
\smallskip \\
\tfrac{\partial^{2}\phi^{2k^{\prime}-1}}
{\partial x^{2j^{\prime}-1}\partial x^{l}}(x)
+\tfrac{\partial^{2}\phi^{2j^{\prime}}} {\partial x^{2k^{\prime}}\partial
x^{l}}(x) = & \lambda_{2j^{\prime}-1,2k^{\prime},l},
\smallskip \\
-\tfrac{\partial^{2}\phi^{2k^{\prime}}} {\partial x^{2j^{\prime}-1}\partial
x^{l}}(x) +\tfrac{\partial^{2}\phi^{2j^{\prime}}} {\partial
x^{2k^{\prime}-1}\partial x^{l}}(x) = &
\lambda_{2j^{\prime}-1,2k^{\prime}-1,l},
\end{array}
\label{syst}
\end{equation}
as follows by taking derivatives with respect to $x^{l}$, $1\leq l\leq2n$, and
evaluating at $x$, in the coefficient of $dx^{j}\wedge dx^{k}$ in the
right-hand side of the first equation in \eqref{Omega_prima}. Furthermore, as a
computation shows, the equations \eqref{compat} are seen to be the
compatibility conditions of the system (\ref{syst}), thus concluding that an
element $j_{x}^{2}\phi\in G_{x}^{21}$ satisfying such a system really exists.

Therefore, a $G_{x}^{2}$-equivariant isomorphism of algebras holds:
\[
\left(  \bar{\delta}_{x}\right) ^{\ast}\colon C^{\infty}((O_{M})_{x}
\times\wedge^{3}T_{x}^{\ast}M)\overset{\cong}{\longrightarrow}C^{\infty}
(J_{x}^{1}O_{M})^{G_{x}^{21}}.
\]

Moreover, as $G_{x}^{2}$ is the semidirect product $G_{x}^{2}=G_{x}
^{21}\rtimes G_{x}^{1}$, taking invariants with respect to $G_{x}^{1}\cong
GL(2n,\mathbb{R})$, we finally obtain an isomorphism
\[
\left(
\bar{\delta}_{x}
\right) ^{\ast}\colon C^{\infty}((O_{M})_{x}
\times\wedge^{3}T_{x}^{\ast}M)^{G_{x}^{1}}
\overset{\cong}{\longrightarrow}
\left( C^{\infty}(J_{x}^{1}O_{M})^{G_{x}^{21}}
\right) ^{G_{x}^{1}}
=C^{\infty}(J_{x}^{1}O_{M})^{G_{x}^{2}}.
\]
Once an element $\Omega_{x}\in(O_{M})_{x}$ has been fixed, the following
injective map is defined:
\[
\begin{array}
[c]{l}
\alpha_{\Omega_{x}}^{1}\colon\wedge^{3}T_{x}^{\ast}M\to(O_{M})_{x}
\times\wedge^{3}T_{x}^{\ast}M,\\
\alpha_{\Omega_{x}}^{1}(\theta_{x})=(\Omega_{x},\theta_{x}),
\end{array}
\]
such that $L_{A}\circ\alpha_{\Omega_{x}}^{1}=\alpha_{A\cdot\Omega_{x}}^{1}
\circ L_{A}$, $\forall A\in G_{x}^{1}$; in particular, if $A\in
Sp(\Omega_{x})$,
we have $L_{A}\circ\alpha_{\Omega_{x}}^{1}=\alpha_{\Omega_{x}}^{1}\circ L_{A}$.
Hence the map $\alpha_{\Omega_{x}}^{1}$ is $Sp(\Omega_{x})$-equivariant. If
$f\in C^{\infty}((O_{M})_{x}\times\wedge ^{3}T_{x}^{\ast}M)^{G_{x}^{1}}$, then
$f\circ\alpha_{\Omega_{x}}^{1}\in
C^{\infty}(\wedge^{3}T_{x}^{\ast}M)^{Sp(\Omega_{x})}$. In fact, if $A\in
Sp(\Omega_{x})$, then for every $\theta_{x}\in\wedge^{3}T_{x}^{\ast}M$ we have
\begin{align*}
\left( f\circ\alpha_{\Omega_{x}}^{1}\right) \left( A\cdot\theta_{x}\right)
&  =f\left(  \Omega_{x},A\cdot\theta_{x}\right) =f\left(  A\cdot\Omega
_{x},A\cdot\theta_{x}\right)  =f\left(  \Omega_{x},\theta_{x}\right) \\ &
=\left(  f\circ\alpha_{\Omega_{x}}^{1}\right)  \left(  \theta_{x}\right) .
\end{align*}
Therefore, the $Sp(\Omega_{x})$-equivariant ring-homomorphism
\[
\left( \alpha_{\Omega_{x}}^{1}\right) ^{\ast}\colon C^{\infty}\left(
(O_{M})_{x}\times\wedge^{3}T_{x}^{\ast}M\right) \to C^{\infty}\left(
\wedge^{3}T_{x}^{\ast}M\right)
\]
induced by $\alpha_{\Omega_{x}}^{1}$ maps $C^{\infty}((O_{M})_{x}
\times \wedge^{3}T_{x}^{\ast}M)^{G_{x}^{1}}$ into
$C^{\infty}(\wedge^{3}T_{x}^{\ast}M)^{Sp(\Omega_{x})}$, and the restriction of
$\left(  \alpha_{\Omega_{x}}^{1}\right) ^{\ast}$ to
$C^{\infty}((O_{M})_{x}\times\wedge^{3}T_{x}^{\ast}M)^{G_{x}^{1}}$ is denoted
by
\[
\left( \alpha_{\Omega_{x}}^{1}\right) ^{\ast\prime}\colon C^{\infty}\left(
(O_{M})_{x}\times\wedge^{3}T_{x}^{\ast}M\right) ^{G_{x}^{1}}\to
C^{\infty}\left(  \wedge^{3}T_{x}^{\ast}M\right) ^{Sp(\Omega_{x})}.
\]

We prove that $\left( \alpha_{\Omega_{x}}^{1}\right) ^{\ast\prime}$ is
injective. Actually, if $\left( \alpha_{\Omega_{x}}^{1}\right)
^{\ast\prime}(f)=0$, then
\[
f\left( \Omega_{x},\theta_{x}\right)  =0,
\qquad\forall\theta_{x}\in\wedge^{3}T_{x}^{\ast}M.
\]
As $f$ is invariant under the action of $G_{x}^{1}$, we also have $f\left(
A\cdot\Omega_{x},A\cdot\theta_{x}\right) =0$, $\forall A\in G_{x}^{1}$,
$\forall\theta_{x}\in\wedge^{3}T_{x}^{\ast}M$; as $G_{x}^{1}$ operates
transitively on $(O_{M})_{x}$, it follows: $f=0$, because given an arbitrary
point $\left( \Omega_{x}^{\prime},\theta_{x}^{\prime}\right)  \in(O_{M})_{x}
\times\wedge^{3}T_{x}^{\ast}M$ there exists $A\in G_{x}^{1}$ such that
$A\cdot\Omega_{x}=\Omega_{x}^{\prime}$ and by taking $\theta_{x}=A^{-1}
\cdot\theta_{x}^{\prime}$, we conclude
$f\left( \Omega_{x}^{\prime},\theta_{x}^{\prime}\right) =0$.

The map $\left( \alpha_{\Omega_{x}}^{1}\right) ^{\ast\prime}$ is also
surjective: For every $g\in C^{\infty}\left( \wedge^{3}T_{x}^{\ast}M\right)
^{Sp(\Omega_{x})}$ we define $f\colon(O_{M})_{x}\times\wedge^{3}T_{x}^{\ast}M
\to\mathbb{R}$ as follows: $f\left( \Omega_{x}^{\prime},\theta_{x}\right)
=g(A^{-1}\cdot\theta_{x})$, $A\in G_{x}^{1}$ being any transformation verifying
$\Omega_{x}^{\prime}=A\cdot\Omega_{x}$. The definition is correct, since if $B$ also verifies the equation
$\Omega_{x}^{\prime}=B\cdot\Omega_{x}$, then $A^{-1}B\in Sp(\Omega_{x})$ and
$g$ being invariant under the action of $Sp(\Omega_{x})$, we have
$g(B^{-1}\cdot\theta_{x})=g((A^{-1}B)^{-1}
\cdot A^{-1}\cdot\theta_{x})=g(A^{-1}\cdot\theta_{x})$. Furthermore, $f$ is
$G_{x}^{1}$-invariant, as $f(A^{\prime}\cdot\Omega_{x}^{\prime},A^{\prime}
\cdot\theta_{x})=g((A^{\prime}A)^{-1}\cdot A^{\prime}\cdot\theta_{x})
=f\left( \Omega_{x}^{\prime},\theta_{x}\right) $, $\forall A^{\prime}\in
G_{x}^{1}$, thus concluding.
\end{proof}
\section{The number of invariants\label{number}}
\subsection{Infinitesimal invariants}
From now onwards, $V$ denotes a real vector space of dimension $2n$ and
$\Omega\in\wedge^{2}V^{\ast}$ denotes a non-degenerate skew-symmetric bilinear
form on $V$.

Let $(v_{i})_{i=1}^{2n}$ be a basis for $V$ with dual basis
$(v^{i})_{i=1}^{2n}$. We define coordinate functions $y_{abc}$, $1\leq
a<b<c\leq2n$, on $\wedge^{3}V^{\ast}$ by setting
\begin{equation}
\theta=\sum\nolimits_{1\leq a<b<c\leq2n}y_{abc}(\theta)\left(  v^{a}\wedge
v^{b}\wedge v^{c}\right) \in\wedge^{3}V^{\ast}. \label{coordinates}
\end{equation}
If $A\in GL(V)$, then for $1\leq a<b<c\leq2n$ we have
\begin{align*}
A\cdot\left(  v^{a}\wedge v^{b}\wedge v^{c}\right) &
=(A^{-1})^{\ast}v^{a}\wedge(A^{-1})^{\ast}v^{b}\wedge(A^{-1})^{\ast}v^{c}\\ &
=\left(  v^{a}\circ A^{-1}\right) \wedge\left(  v^{b}\circ A^{-1}\right)
\wedge\left(  v^{c}\circ A^{-1}\right) \\
&  =\left(  \lambda_{h}^{a}v^{h}\right)  \wedge
\left(  \lambda_{i}^{b}v^{i}\right)
\wedge\left(  \lambda_{j}^{c}v^{j}\right) \\
&  =\sum\nolimits_{1\leq h<i<j\leq2n}\left\vert
\begin{array}
[c]{ccc}
\lambda_{ah} & \lambda_{bh} & \lambda_{ch}\\
\lambda_{ai} & \lambda_{bi} & \lambda_{ci}\\
\lambda_{aj} & \lambda_{bj} & \lambda_{cj}
\end{array}
\right\vert v^{h}\wedge v^{i}\wedge v^{j},
\end{align*}
where $(\lambda_{ij})_{i,j=1}^{2n}$ is the matrix of $(A^{-1})^{T}$ in the
basis $(v_{i})_{i=1}^{2n}$ and the superscript $T$ means transpose. In what
follows we assume $\Omega=\sum_{i=1}^{n}v^{i}\wedge v^{n+i}$.

A function $I\colon\wedge^{3}V^{\ast}\to\mathbb{R}$ is $Sp(\Omega )$-invariant
if $I\left(  \Lambda\cdot\theta\right)  =I(\theta)$,
$\forall\theta\in\wedge^{3}V^{\ast}$, $\forall\Lambda\in Sp(\Omega) $.
\begin{lemma}
\label{lemma_bis}A smooth function $I\colon\wedge^{3}V^{\ast}
\to \mathbb{R}$ is $Sp(\Omega)$-invariant if and only if $I$
is a first integral of the distribution spanned by the following vector
fields:
\begin{equation}
\begin{array}
[c]{rc} U^{\ast}= & \sum\nolimits_{1\leq h<i<j\leq2n}\left(
\sum\nolimits_{1\leq a<b<c\leq2n}U_{hij}^{abc}y_{abc}\right)
\tfrac{\partial}{\partial y_{hij}},
\smallskip\\
U= & \multicolumn{1}{l}{(u_{ij})_{i,j=1}^{2n}\in\mathfrak{sp}(2n,\mathbb{R}),}
\end{array}
\label{U^ast}
\end{equation}
where the functions $U_{hij}^{abc}$ are given by the formulas
\[
U_{hij}^{abc}=-\left\vert
\begin{array}
[c]{ccc} u_{ha} & \delta_{hb} & \delta_{hc}\\ u_{ia} & \delta_{ib} &
\delta_{ic}\\ u_{ja} & \delta_{jb} & \delta_{jc}
\end{array}
\right\vert -\left\vert
\begin{array}
[c]{ccc}
\delta_{ha} & u_{hb} & \delta_{hc}\\
\delta_{ia} & u_{ib} & \delta_{ic}\\
\delta_{ja} & u_{jb} & \delta_{jc}
\end{array}
\right\vert -\left\vert
\begin{array}
[c]{ccc}
\delta_{ha} & \delta_{hb} & u_{hc}\\
\delta_{ia} & \delta_{ib} & u_{ic}\\
\delta_{ja} & \delta_{jb} & u_{jc}
\end{array}
\right\vert .
\]
\end{lemma}
\begin{proof}
If $I$ is invariant, then, in particular, we have $I\left(
\exp(tU)\cdot\theta\right)  =I(\theta)$, $\forall t\in\mathbb{R}$,
$U=(u_{ij})_{i,j=1}^{2n}\in\mathfrak{sp}(\Omega)$. If
$\Lambda(t)=\exp(-tU^{T})$, then
\[
I\Bigl(
\sum\nolimits_{\substack{1\leq a<b<c\leq2n\\1\leq h<i<j\leq2n}}~y_{abc}
\left\vert
\begin{array}
[c]{ccc}
\lambda_{ah}(t) & \lambda_{bh}(t) & \lambda_{ch}(t)\\
\lambda_{ai}(t) & \lambda_{bi}(t) & \lambda_{ci}(t)\\
\lambda_{aj}(t) & \lambda_{bj}(t) & \lambda_{cj}(t)
\end{array}
\right\vert v^{h}\wedge v^{i}\wedge v^{j}
\Bigr)
=I(\theta),
\]
and taking derivatives at $t=0$, it follows:
\[
0=\sum\nolimits_{1\leq a<b<c\leq6,1\leq h<i<j\leq6}U_{hij}^{abc}y_{abc}
\tfrac{\partial I}{\partial y_{hij}}(\theta),
\]
$U_{hij}^{abc}$ being as in the statement. The converse follows from the fact
that the symplectic group is connected and hence, every symplectic
transformation is a product of exponentials of matrices in the symplectic
algebra.
\end{proof}

\begin{theorem}
\label{Th2}The distribution $\mathcal{D}\subset
T\mathfrak{(}\wedge^{3}V^{\ast})$
whose fibre $\mathcal{D}_{\theta}$\ over $\theta\in\wedge ^{3}V^{\ast}$ is the
subspace $(U^{\ast})_{\theta}$, $U\in\mathfrak{sp}(2n,\mathbb{R})$, is
involutive and of locally constant rank on a dense open subset
$\mathcal{O}\subset\wedge^{3}V^{\ast}$.

The number $N_{2n}$ of functionally independent $Sp(\Omega)$-invariant
functions defined on $\mathcal{O}$ is equal to
$N_{2n}=\tbinom{2n}{3}-\operatorname*{rank}\mathcal{D}|_{\mathcal{O}}$.
\end{theorem}
\begin{proof}
Every pair of vector fields $U^{\prime\ast},U^{\prime\prime\ast}$ belonging to
$\mathcal{D}$ on an open subset $O\subseteq\wedge^{3}V^{\ast}$ can be written
as $U^{\prime\ast}=\sum_{h=1}^{n(2n+1)}f^{h}(U_{h})^{\ast}$, $U^{\prime
\prime\ast}=\sum_{i=1}^{n(2n+1)}g^{i}(U_{i})^{\ast}$, with $f^{h},g^{i}\in
C^{\infty}(O)$, where $(U_{1},\ldots,U_{n(2n+1)})$ is a basis of
$\mathfrak{sp}(2n,\mathbb{R)}$. As
$[U_{h}^{\ast},U_{i}^{\ast}]=-[U_{h},U_{i}]^{\ast}$, it follows that
$[U^{\prime\ast},U^{\prime\prime\ast}]$ can be written as a linear combination
of $(U_{h})^{\ast}$, $(U_{i})^{\ast}$, and
$[U_{h}^{\ast},U_{i}^{\ast}]=-[U_{h},U_{i}]^{\ast}=-c_{hi}^{j}(U_{j})^{\ast}$,
where $c_{hi}^{j}$ are the structure constants of
$\mathfrak{sp}(2n,\mathbb{R})$ on this basis. This proves that $\mathcal{D}$ is
involutive. Moreover, we first recall that the dimension of the vector spaces
$\{\mathcal{D}_{\theta}:\theta\in\wedge^{3}V^{\ast}\}$ is uniformly bounded by
$\dim(\wedge^{3}V^{\ast})=\tbinom{2n}{3}$. Let
$\mathcal{O}\subset\wedge^{3}V^{\ast}$ be the subset defined as follows: A
point $\theta\in\wedge^{3}V^{\ast}$ belongs to $\mathcal{O}$ if and only if
$\theta$ admits an open neighbourhood $N$ such that $d=\dim\mathcal{D}_{\theta}
=\max_{\theta^{\prime}\in N}\left( \dim\mathcal{D}_{\theta^{\prime}}\right) $.
We claim that $\mathcal{O}$ is an open subset. Actually, there exists an open
neighbourhood $N^{\prime}$ of $\theta$ such that the dimension of the fibres of
$\mathcal{D}$ over the points $\theta^{\prime}\in N^{\prime}$ is at least $d$,
as if $\left.  (X_{i})^{\ast}\right\vert _{\theta}$, $1\leq i\leq d$, is a
basis for $\mathcal{D}$ at $\xi$, for certain
$X_{i}\in\mathfrak{sp}(2n,\mathbb{R})$, $1\leq i\leq d$, then the vector fields
$(X_{i})^{\ast}$, $1\leq i\leq d$, are linearly independent at each point of a
neighbourhood of $\theta$. From the definition of $\mathcal{O}$ we thus
conclude that if $\theta\in\mathcal{O}$, then we have
$\dim\mathcal{D}_{\theta^{\prime}}=d$ for every $\theta^{\prime}\in N\cap
N^{\prime}$; hence $N\cap N^{\prime}\subseteq\mathcal{O}$. The same argument
proves that the rank of $\mathcal{D}$ is locally constant over $\mathcal{O}$.
Next, we prove that $\mathcal{O}$ is dense. Let $N$ be an open neighbourhood of
an arbitrarily chosen point $\theta\in\wedge^{3}V^{\ast}$ and let
$\theta^{\prime}$ be a point in $N$ such that the rank of $\mathcal{D}|_{N}$
takes its greatest value at $\theta^{\prime}$. By proceeding as above, we
deduce that $\theta^{\prime}$ belongs to $\mathcal{O}$. Finally, the formula
for the number of invariants in the statement now follows from Frobenius'
theorem.
\end{proof}
\begin{remark}
\label{remark1}
We have $N_{2}=0$, as $\wedge^{3}V^{\ast}=\{0\}$ if $\dim V=2$, and $N_{4}=0$,
as $Sp(2n)$ acts transitively on $\wedge^{3}V^{\ast}
\backslash\{0\}$ if $\dim V=2n=4$. Furthermore, as a consequence
of the results obtained in \cite{JL}, it follows that the generic rank of
$\mathcal{D}$ for $\dim V=2n=6$ is $18$; hence $N_{6}=2$.
\end{remark}
\subsection{$N_{2n}$ computed}
\begin{theorem}
\label{Th3}We have
\[
N_{2n}=\left\{
\begin{array}
[c]{rl} 0, & 1\leq n\leq2,\\ 2, & n=3,\\
\frac{n(4n^{2}-12n-1)}{3}, & n\geq4.
\end{array}
\right.
\]
\end{theorem}
\begin{proof}
The formula in the statement for $1\leq n\leq3$ follows from Remark
\ref{remark1}. Hence we can assume $n\geq4$. For every $3$-covector
$\theta \in\wedge^{3}V^{\ast}$, let us define
\[
\begin{array}
[c]{lll}
\mu_{\theta}\colon Sp(\Omega)\to\wedge^{3}V^{\ast},
& \mu_{\theta}(U)=U\cdot\theta, & \forall U\in Sp(\Omega).
\end{array}
\]
It suffices to prove that there exists a dense open subset
$\mathcal{O}^{\prime}\subset\wedge^{3}V^{\ast}$ such that $\mu_{\theta}$ is an
immersion if $\theta\in\mathcal{O}^{\prime}$, since $\operatorname*{im}\left[
(\mu_{\theta})_{\ast}\colon T_{I}Sp(\Omega)\to T_{\theta}(\wedge
^{3}V^{\ast})\right] =\mathcal{D}_{\theta}$.

By expanding on \eqref{U^ast} it follows:
\begin{equation}
\begin{array}
[c]{l} -U^{\ast}=\\
\sum\limits_{1\leq a<b<c\leq2n}\;\sum\limits_{1\leq d<b}u_{da}y_{abc}Y_{dbc}
+\sum\limits_{1\leq a<b<c\leq2n}\;\sum\limits_{c<d\leq2n}u_{da}y_{abc}Y_{bcd}
\medskip\\
-\sum\limits_{1\leq a<b<c\leq2n}\;\sum\limits_{b<d<c}u_{da}y_{abc}Y_{bdc}
+\sum\limits_{1\leq a<b<c\leq2n}\;\sum\limits_{a<d<c}u_{db}y_{abc}Y_{adc}
\medskip\\
-\sum\limits_{1\leq a<b<c\leq2n}\;\sum\limits_{c<d\leq2n}u_{db}y_{abc}Y_{acd}
-\sum\limits_{1\leq a<b<c\leq2n}\;\sum\limits_{1\leq d<a}u_{db}y_{abc}Y_{dac}
\medskip\\
+\sum\limits_{1\leq a<b<c\leq2n}\;\sum\limits_{b<d\leq2n}u_{dc}y_{abc}Y_{abd}
+\sum\limits_{1\leq a<b<c\leq2n}\;\sum\limits_{1\leq d<a}u_{dc}y_{abc}Y_{dab}
\medskip\\
-\sum\limits_{1\leq a<b<c\leq2n}\;\sum\limits_{a<d<b}u_{dc}y_{abc}Y_{adb},
\end{array}
\label{U^ast_2}
\end{equation}
where $Y_{hij}=\tfrac{\partial}{\partial y_{hij}}$, $1\leq h<i<j\leq2n$. Given
indices $1\leq\alpha<\beta<\gamma \leq2n$, the coefficient of
$Y_{\alpha\beta\gamma}$ in \eqref{U^ast_2} is
\[
\begin{array}
[c]{rl} C_{\alpha\beta\gamma}= & \sum\nolimits_{a=1}^{\beta-1}u_{\alpha
a}y_{a\beta\gamma} +\sum\nolimits_{a=1}^{\beta-1}u_{\gamma a}y_{a\alpha\beta
}-\sum\nolimits_{a=1}^{\alpha-1}u_{\beta a}y_{a\alpha\gamma}\medskip\\ &
+\sum\nolimits_{a=\alpha+1}^{\gamma-1}u_{\beta a}y_{\alpha a\gamma}
-\sum\nolimits_{a=\alpha+1}^{\beta-1}u_{\gamma a}y_{\alpha a\beta}
-\sum\nolimits_{a=\beta+1}^{\gamma-1}u_{\alpha a}y_{\beta a\gamma}
\medskip\\
& +\sum\nolimits_{a=\beta+1}^{2n}u_{\gamma a}y_{\alpha\beta a} +\sum
_{a=\gamma+1}^{2n}u_{\alpha a}y_{\beta\gamma a} -\sum_{a=\gamma+1}^{2n}u_{\beta
a}y_{\alpha\gamma a},
\medskip\\
& \multicolumn{1}{r}{1\leq\alpha<\beta<\gamma\leq2n.}
\end{array}
\]
As the matrix $U=(u_{ij})_{i,j=1}^{2n}$ is symplectic, the following symmetries
hold:
\[
\begin{array}
[c]{rrrl} u_{j,n+i}=u_{i,n+j}, & u_{n+j,i}=u_{n+i,j}, & u_{n+j,n+i}=-u_{i,j}, &
1\leq i<j\leq n.
\end{array}
\]

A vector $U_{\theta}^{\ast}$ belongs to $\ker(\mu_{\theta})_{\ast}$ if and only
if $C_{\alpha\beta\gamma}=0$ for every system of indices $1\leq
\alpha<\beta<\gamma\leq2n$. We thus obtain a homogeneous linear system
$S_{2n}$\ of $\tbinom{2n}{3}$\ linear equations in the $n(2n+1)$ unknowns
$u_{ij}$, $i,j=1,\ldots,n$; $u_{i,n+j}$, $u_{n+i,j}$, $1\leq i\leq j\leq n$,
and we have $\tbinom{2n}{3}>n(2n+1)$ for every $n\geq4$. Evaluating $S_{2n}$ at
the $3$-covector $\theta^{0}$ of coordinates $y_{abc}(\theta^{0})=a+b+c$,
$1\leq a<b<c\leq2n$, as a numerical calculation shows, the only solution to
$S_{2n}(\theta^{0})$ is given by $u_{ij}=0$, $i,j=1,\ldots,n$. We can thus
conclude by simply applying the formula \eqref{Th2} for $N_{2n}$\ in Theorem
\ref{Th2}.
\end{proof}


\begin{thebibliography}{9}
\bibitem [1]{JL} J. Mu\~{n}oz Masqu\'e, L. M. Pozo Coronado, \emph{A new look at
the classification of the tri-covectors of a }$6$\emph{-dimensional symplectic
space}, Linear and Multilinear Algebra, {\bf 67} (5), 2019, pp. 939--952.
DOI:doi.org/10.1080/03081087.2018.1440517.
\end{thebibliography}
\end{document}